\documentclass[11pt]{article}
\usepackage[margin=0.9in,headheight=14pt]{geometry}
\usepackage[T1]{fontenc}
\usepackage{lmodern,microtype,amsmath,amssymb,amsthm,mathtools,nicefrac}
\usepackage{xcolor,xspace,tikz,booktabs,needspace,fancyhdr}
\usepackage{caption}
\usepackage[authoryear,round]{natbib}
\usepackage[framemethod=tikz]{mdframed}
\definecolor{ink}{HTML}{243B40}
\definecolor{jade}{HTML}{297F79}
\definecolor{ochre}{HTML}{B8792D}
\usepackage[colorlinks=true,linkcolor=jade!80!black,citecolor=jade!80!black,urlcolor=jade!80!black]{hyperref}
\usetikzlibrary{arrows.meta,calc}
\newcommand{\algname}[1]{{\color[gray]{0.35}\small\sf#1}\xspace}
\newcommand{\eqdef}{\coloneqq}
\newcommand{\X}{\mathcal X}
\newcommand{\B}{\mathbb B}
\newcommand{\ip}[2]{\langle #1,#2\rangle}
\newcommand{\norm}[1]{\lVert #1\rVert}
\newcommand{\R}{\mathbb R}
\newcommand{\argmin}{\operatorname*{arg\,min}}
\newtheorem{assumption}{Assumption}
\newmdtheoremenv[backgroundcolor=black!3,linewidth=0pt,innertopmargin=6pt,innerbottommargin=6pt,nobreak=true,skipabove=10pt,skipbelow=9pt]{theorem}{Theorem}
\newmdtheoremenv[backgroundcolor=black!3,linewidth=0pt,innertopmargin=6pt,innerbottommargin=6pt,nobreak=true,skipabove=10pt,skipbelow=9pt]{proposition}[theorem]{Proposition}
\newmdtheoremenv[backgroundcolor=black!3,linewidth=0pt,innertopmargin=6pt,innerbottommargin=6pt,nobreak=true,skipabove=10pt,skipbelow=9pt]{corollary}[theorem]{Corollary}
\newmdtheoremenv[backgroundcolor=black!3,linewidth=0pt,innertopmargin=6pt,innerbottommargin=6pt,nobreak=true,skipabove=10pt,skipbelow=9pt]{lemma}[theorem]{Lemma}
\title{\bfseries\color{ink}Frank--Wolfe is Approximate\\Broximal Gradient Descent}
\author{Peter Richt\'arik\qquad Zhirayr Tovmasyan\\[3pt]\normalsize King Abdullah University of Science and Technology}
\date{\small November 1, 2026}
\begin{document}
\maketitle
\begin{abstract}
We show that the Frank--Wolfe method is an approximation of Broximal Gradient Descent (BroxGD), and hence of the ball-proximal method. For constrained minimization, the ball-proximal method minimizes the objective over the intersection of the feasible set with a ball around the current iterate. Linearizing the objective gives BroxGD, but the intersection can remain difficult to handle even when linear minimization over the feasible set alone is cheap. We therefore replace the intersection by a smaller copy of the feasible set, obtained by shrinking it toward the current iterate so that it lies inside the ball. Minimizing the linear model over this smaller set gives exactly Frank--Wolfe, with the shrinking factor as its step size.

We examine what this approximation preserves and what it can lose. For smooth convex objectives on compact convex sets, Frank--Wolfe retains a linear-model decrease guarantee that gives it and BroxGD a common convergence proof, recovering the classical $O(1/k)$ bound. Among compact convex regions inside the ball chosen before the gradient is known, the smaller copy is the smallest that preserves the same linear-model decrease guarantee. Nevertheless, Frank--Wolfe can achieve much less linear-model decrease than BroxGD at the same radius. A simplex example illustrates the sparsity gained by restricting the region, while a family of triangles shows how Frank--Wolfe can lose the distance decrease guaranteed by BroxGD. We also give a computable test for recovering the distance guarantee and a refinement procedure that terminates after finitely many oracle calls on polytopes.
\end{abstract}

\section{Introduction}\label{sec:start}
We consider the constrained minimization problem
\begin{equation}\label{eq:problem}
 \min_{x\in\X} f(x),
\end{equation}
where $\X\subset\R^n$ is a nonempty compact convex set and $f$ is convex and differentiable with an $L$-Lipschitz gradient on $\X$. Our starting point is the ball-proximal method (\algname{BPM}) of \citet{brox}. For a proper closed function $F$, a starting point $u_0\in\operatorname{dom}F$, and radii $t_k>0$, it sets
\begin{equation}\label{eq:bpm}
 u_{k+1}\in\argmin_{z\in\B(u_k,t_k)}F(z);
\end{equation}
for a fixed radius, the map sending the current point to the set of minimizers in \eqref{eq:bpm} is the broximal operator of \citet{brox}. For a proper closed convex $F$ with a minimizer and a constant radius, the objective gap after $K$ iterations decays like $\exp(-cK^2)$, with $c$ determined by the radius and the initial distance to the solution set, and a minimizer is reached in finitely many iterations (\citealp[Theorem~8.1]{brox}; \citealp[Theorem~3.4]{sharper}). Neither smoothness nor strong convexity is required. Such guarantees cannot hold uniformly over general convex objectives in the standard first-order oracle model \citep{nemirovski}; \algname{BPM} obtains them through a stronger oracle, the exact minimization of the objective over a ball.

That oracle can make the steps of \algname{BPM} impractical to compute: if the ball contains a minimizer of $F$, every solution of \eqref{eq:bpm} is a global minimizer \citep[Theorem~8.1]{brox}, so the step is in general no easier than the original problem. \algname{BPM} is therefore a conceptual method, and implementable algorithms are obtained by approximating its step (\citealp[Section~1.5]{brox}; \citealp{warsawtalk}). In this paper we follow this route for problem \eqref{eq:problem}, and arrive at the Frank--Wolfe method.

\begin{figure}[t]
\begin{center}
\begin{tikzpicture}[scale=0.96,>=Latex,font=\small]
% Convex octagon with eight exposed vertices; D joins v0 and v4.
% x=(.9,1.4), gamma=.32, g=(-1,.55), t=gamma*sqrt(31.025).
% b=x-t*g/||g|| is inside X; s=v3 is the unique global minimizer.
\def\covergamma{0.32}
\pgfmathsetmacro{\coverradius}{\covergamma*sqrt(31.025)}
\coordinate (x) at (.9,1.4);
\coordinate (v0) at (-.35,.65);
\coordinate (v1) at (.15,-.5);
\coordinate (v2) at (2.55,-.7);
\coordinate (v3) at (4.85,.1);
\coordinate (v4) at (5.1,1.8);
\coordinate (v5) at (3.9,3.15);
\coordinate (v6) at (1.4,3.35);
\coordinate (v7) at (-.2,2.1);
\def\coverpoly{(v0)--(v1)--(v2)--(v3)--(v4)--(v5)--(v6)--(v7)--cycle}
\foreach \i in {0,...,7}{\coordinate (h\i) at ($(x)!\covergamma!(v\i)$);}
\coordinate (s) at (v3);
\coordinate (y) at (h3);
\coordinate (b) at ({.9+\coverradius/sqrt(1+.55^2)},
                         {1.4-.55*\coverradius/sqrt(1+.55^2)});
% The tinted intersection and its ball boundary are drawn separately.
\fill[ink!4] \coverpoly;
\draw[jade!40,densely dashed,line width=.55pt] (x) circle (\coverradius);
\begin{scope}
\clip \coverpoly;
\fill[jade!12] (x) circle (\coverradius);
\draw[jade!70,line width=.7pt] (x) circle (\coverradius);
\end{scope}
\draw[ink!55,line width=.75pt] \coverpoly;
% Two faint rays make the homothety visible without obscuring the steps.
\foreach \i in {5,7}{\draw[ink!18,densely dashed,line width=.45pt] (x)--(v\i);}
\filldraw[fill=ochre!23,draw=ochre,line width=.9pt]
 (h0)--(h1)--(h2)--(h3)--(h4)--(h5)--(h6)--(h7)--cycle;
\foreach \i in {0,...,7}{
 \fill[ink!55] (v\i) circle (1.25pt);
 \fill[ochre!85!black] (h\i) circle (.8pt);
}
\draw[->,ink!40,densely dashed,line width=.65pt] (x)--(s);
\draw[->,jade,line width=1.1pt] (x)--(b);
\draw[->,ochre,line width=1.35pt] (x)--(y);
\fill[ink] (x) circle (2pt) node[above=4pt,text=ink] {$x$};
\fill[ochre] (y) circle (2.4pt) node[above right=1pt,text=ochre!85!black] {$y$};
\fill[jade] (b) circle (2.4pt) node[below right=1pt,text=jade] {$b$};
\fill[ink] (s) circle (2pt) node[below right=2pt,text=ink] {$s$};
\node[text=ink!75] at (4.2,2.4) {$\X$};
\node[text=jade] at (.65,2.65) {$\mathcal C_t(x)$};
\draw[ochre!75,line width=.5pt] (1.86,1.96)--(2.5,2.04);
\node[anchor=west,text=ochre!85!black] at (2.5,2.04) {$\mathcal H_\gamma(x)$};
\node[text=ink!70] at (2.1,-1.06)
 {$\mathcal H_\gamma(x)\subseteq\mathcal C_t(x)=\X\cap\B(x,t)$};
% The right-hand derivation uses the same colors as the two feasible regions.
\node[anchor=west,text=ink] at (6,3.08)
 {\algname{BPM}: $\displaystyle\min_{z\in\B(x,t)}\{f(z)+\iota_{\X}(z)\}$};
\draw[->,ink!35,line width=.7pt] (6.25,2.66)--(6.25,1.98);
\node[anchor=west,font=\small\sffamily,text=ink!70] at (6.6,2.32)
 {Linearize $f$; keep $\iota_{\X}$.};
\node[anchor=west,text=jade] at (6,1.57)
 {\algname{BroxGD}: $\displaystyle b\in\argmin_{z\in\mathcal C_t(x)}\langle g,z-x\rangle$};
\draw[->,ochre!80,line width=.9pt] (6.25,1.15)--(6.25,.47);
\node[anchor=west,font=\small\sffamily,text=ochre!85!black] at (6.6,.81)
 {Replace by contracted $\X$.};
\node[anchor=west,text=ochre!85!black] at (6,.06)
 {\algname{FW}: $\displaystyle y\in\argmin_{z\in\mathcal H_\gamma(x)}\langle g,z-x\rangle$};
\node[anchor=west,text=ink!70] at (6,-.82)
 {$y=x+\gamma(s-x),\qquad t=\gamma D.$};
\end{tikzpicture}
\end{center}
\caption{The two approximations of a \algname{BPM} step at a point $x\in\X$ with gradient $g$. Left: the feasible set $\X$, the ball intersection $\mathcal C_t(x)=\X\cap\B(x,t)$ of \algname{BroxGD}, and the contracted set $\mathcal H_\gamma(x)=x+\gamma(\X-x)$, with $t=\gamma D$. The point $s$ is the global linear minimizer over $\X$, $b$ is the exact \algname{BroxGD} step \eqref{eq:brox}, and $y=x+\gamma(s-x)$ is the Frank--Wolfe step \eqref{eq:fw-step}, which minimizes the same linear objective over $\mathcal H_\gamma(x)$. Right: the corresponding chain of subproblems.}
\label{fig:regions}
\end{figure}
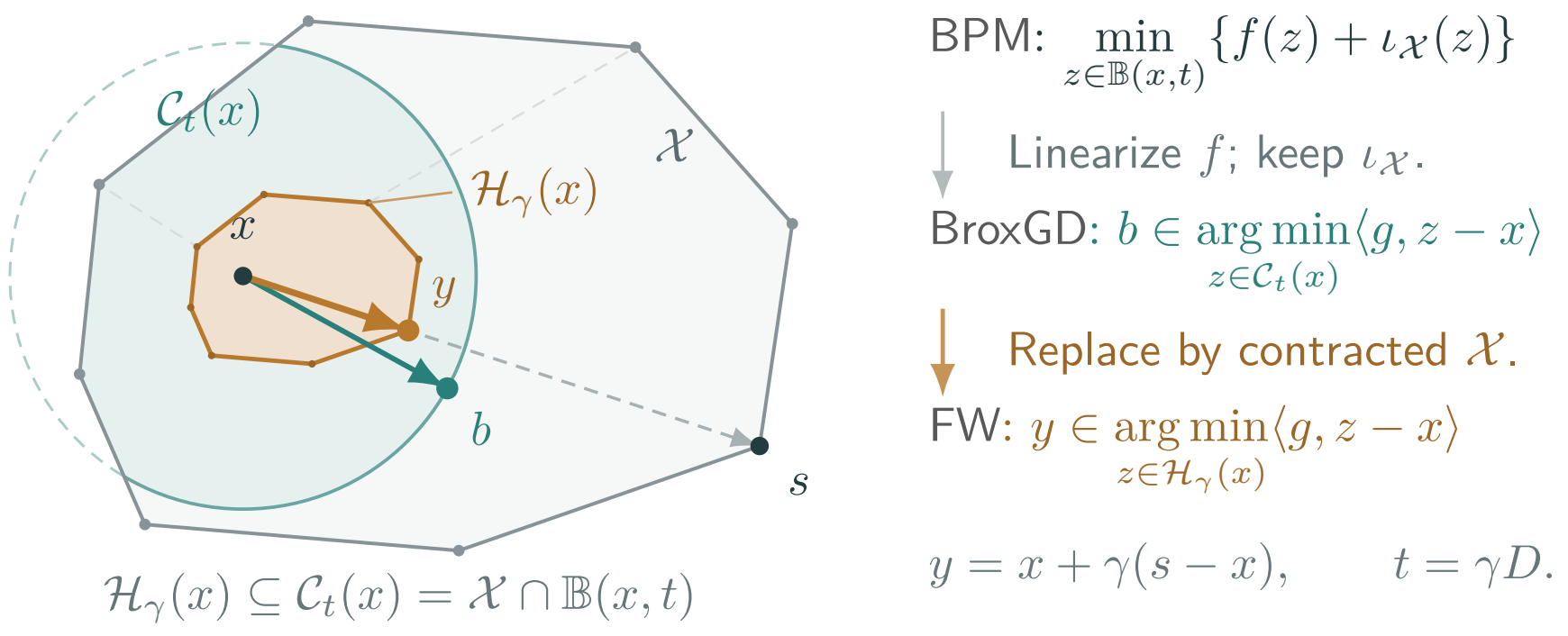

The natural first approximation replaces the smooth part of the objective by its linear model and keeps the rest. For problem \eqref{eq:problem}, the constraint enters the objective through the indicator function $\iota_{\X}$ of $\X$ ($0$ on $\X$, $+\infty$ outside), so that \eqref{eq:problem} is the minimization of $F=f+\iota_{\X}$. With $g=\nabla f(x)$, the step at the current point $x$ becomes
\begin{equation}\label{eq:brox}
 b\in\argmin_{z\in\mathcal C_t(x)}\ip g{z-x},\qquad
 \mathcal C_t(x)\eqdef\X\cap\B(x,t).
\end{equation}
This is Broximal Gradient Descent (\algname{BroxGD}) of \citet{broxgd}, the forward--backward method whose backward step applies the broximal operator to the linearized composite model.

The ball bounds the error of the linear model on the region where it is used, ensures that \eqref{eq:brox} has a solution even when $\X$ is unbounded, and yields convergence rates that do not depend on the diameter of $\X$ \citep[Table~1]{broxgd}. However, the intersection in the \algname{BroxGD} subproblem can make the step difficult to compute, even when linear minimization over $\X$ and over $\B(x,t)$ is tractable separately \citep[Appendix~B.5]{broxgd}. For example, on a nuclear-norm ball the global oracle uses only a leading singular-vector pair of the gradient \citep[Table~1]{jaggi}, whereas a direct computation of the local step can involve a full singular value decomposition (Appendix~\ref{sec:nuclear-example}).

This difficulty invites one more approximation, now of the feasible region. The ball stays, since it controls the movement and the error of the linear model. Inside it we seek a region over which linear minimization is no harder than over $\X$ itself; we call the oracle that minimizes linear functions over $\X$ the global oracle. Since the global oracle also solves linear minimization over every translate and scaling of $\X$, the natural candidate is the copy of $\X$ scaled toward $x$ by a factor $\gamma\in[0,1]$, the contracted set $\mathcal H_\gamma(x)\eqdef x+\gamma(\X-x)$. It lies in $\X$ by convexity and in $\B(x,t)$ whenever $\gamma D\le t$, where $D$ bounds the diameter of $\X$; we use the correspondence $t=\gamma D$. Linear minimization over $\mathcal H_\gamma(x)$ is solved exactly by scaling the global oracle output $s\in\argmin_{z\in\X}\ip g z$:
\begin{equation}\label{eq:fw-step}
 y=x+\gamma(s-x).
\end{equation}
This is the Frank--Wolfe step \citep{frankwolfe,jaggi}, with step fraction $\gamma=t/D$ equal to the contraction ratio. Figure~\ref{fig:regions} shows the three nested regions. In the classical derivation of Frank--Wolfe, the linear model is minimized over the whole feasible set, which gives the direction $s-x$, and the step size along this direction is chosen afterwards, by a schedule or a line search, because the model cannot be trusted far from $x$ (\citealp[Section~2.2]{bertsekas}; \citealp{jaggi}). The region where the model is trusted thus enters only through the step size, after the direction has been fixed. For a prescribed radius $t=\gamma D$, our derivation fixes the local region first and obtains the direction and step fraction from it. Frank--Wolfe thus arises as an approximation of \algname{BroxGD}, and through it of the ball-proximal method.

It remains to ask what such a coarse approximation retains. A natural measure of a step's quality is the ratio of the decrease of the linear model that the step achieves, $\ip g{x-y}$, to the best decrease available within the ball, $\max_{z\in\mathcal C_t(x)}\ip g{x-z}$. For the Frank--Wolfe step this ratio can be as small as $\gamma$, and Section~\ref{sec:discover} shows that a sequence attaining exactly this ratio at every step need not converge to a minimizer. What the step retains is the inequality
\begin{equation}\label{eq:comparison}
 \ip g{y-x}\le\gamma\ip g{z-x}\qquad(z\in\X),
\end{equation}
which we call the contracted comparison: in the linear model, $y$ is at least as good as the contraction $x+\gamma(z-x)$ of every feasible point. \algname{BroxGD} satisfies \eqref{eq:comparison} because its feasible set contains every such contraction. The convergence proof needs this comparison at the unknown minimizer $x_\star$, which is why we retain it for every feasible point. Convexity and smoothness then give a common recurrence for the objective gap, and suitable step sizes yield the classical $O(1/k)$ bound for both methods (Section~\ref{sec:convergence}). What the approximation gains is sparsity: the simplex example in Section~\ref{sec:discover} gives a dense \algname{BroxGD} step from a vertex, while the Frank--Wolfe step using a vertex-returning oracle has only two nonzero coordinates, the sparsity property of Frank--Wolfe \citep{jaggi}. What is lost is the geometry of the \algname{BroxGD} step: on a triangle, \algname{BroxGD} reaches the minimizer in one step while the Frank--Wolfe step moves away from it (Section~\ref{sec:geometry}). Section~\ref{sec:extension} and Appendix~\ref{sec:geometric-appendix} give a computable accuracy condition under which a refined local solve recovers the distance guarantee.

Section~\ref{sec:discover} proves the facts used above and compares the two steps on the simplex; Section~\ref{sec:convergence} extends the shared convergence proof to inexact oracles, and Section~\ref{sec:radius} treats adaptive step sizes within the same framework.

The contracted set has appeared before, with a different motivation: \citet{nesterov2018} introduced it as a trust region for a composite Frank--Wolfe variant, and \citet{doikov} built affine-invariant higher-order methods on it; we return to this connection in Section~\ref{sec:perspective}. The algorithms and rates are classical \citep{frankwolfe,dunn,guelat,jaggi,lacoste}, and local linear subproblems of the form \eqref{eq:brox} go back to \citet{ferris} and underlie the local oracle of \citet{garber}. The contribution of this paper is another view of Frank--Wolfe, as an approximation of the ball-proximal method through \algname{BroxGD}. Comparing the two steps at the same radius explains what the approximation preserves and what it loses.

\section{The contracted set and what it preserves}\label{sec:discover}
Throughout, $\ip\cdot\cdot$ and $\norm\cdot$ denote the Euclidean inner product and norm, $\B(x,t)$ is the closed ball of radius $t$ centered at $x$, and $x_\star$ is a minimizer of \eqref{eq:problem}.
\begin{assumption}[Compact convex feasible set]\label{ass:set}
The set $\X\subset\R^n$ is nonempty, compact, and convex, and a constant $D>0$ satisfies
\begin{equation}\label{eq:diameter}
 \norm{u-v}\le D\qquad(u,v\in\X).
\end{equation}
\end{assumption}

Throughout this section $x\in\X$ is the current point, $g\in\R^n$ is an arbitrary vector (in the algorithm, $g=\nabla f(x)$), $\gamma\in[0,1]$, and $t\ge0$. The comparison involves the following quantities. The global oracle returns
\begin{equation}\label{eq:oracle}
 s\in\argmin_{z\in\X}\ip g z,\qquad G\eqdef\ip g{x-s},
\end{equation}
where $G$ is the largest decrease of the linear model over $\X$, the Frank--Wolfe gap. The exact local step \eqref{eq:brox} attains the largest decrease over the ball intersection,
\begin{equation}\label{eq:local-value}
 M_t(x)\eqdef\max_{z\in\mathcal C_t(x)}\ip g{x-z}=\ip g{x-b}.
\end{equation}
Both are nonnegative, since $x$ is feasible for both problems.

\begin{lemma}[Geometry of the contraction]\label{lem:contraction}
Let Assumption~\ref{ass:set} hold, and fix $x\in\X$ and $\gamma\in[0,1]$. The contracted set
\begin{equation}\label{eq:homothety}
 \mathcal H_\gamma(x)=x+\gamma(\X-x)=\{x+\gamma(z-x):z\in\X\}
\end{equation}
satisfies
\begin{equation}\label{eq:inclusion}
 \mathcal H_\gamma(x)\subseteq\mathcal C_{\gamma D}(x)\subseteq\X.
\end{equation}
For every $g\in\R^n$, the point $y=x+\gamma(s-x)$ with $s$ from \eqref{eq:oracle} minimizes $\ip g{u-x}$ over $u\in\mathcal H_\gamma(x)$, and
\begin{equation}\label{eq:contracted-value}
 \min_{u\in\mathcal H_\gamma(x)}\ip g{u-x}=\gamma\min_{z\in\X}\ip g{z-x}=-\gamma G.
\end{equation}
\end{lemma}
\begin{proof}
For $z\in\X$, the point $x+\gamma(z-x)$ lies in $\X$ by convexity, and
\[
 \norm{x+\gamma(z-x)-x}=\gamma\norm{z-x}\overset{\eqref{eq:diameter}}{\le}\gamma D.
\]
This proves \eqref{eq:inclusion}.

Every $u\in\mathcal H_\gamma(x)$ has the form $u=x+\gamma(z-x)$ for some $z\in\X$.
Since $\gamma\ge0$ and $s$ solves \eqref{eq:oracle},
\[
 \ip g{u-x}=\gamma\ip g{z-x}
 \overset{\eqref{eq:oracle}}{\ge}\gamma\ip g{s-x}=\ip g{y-x}=-\gamma G.
\]
Thus $y$ minimizes the linear model over $\mathcal H_\gamma(x)$,
and \eqref{eq:contracted-value} follows.
\end{proof}

The next result compares the linear-model decreases of the Frank--Wolfe step and the exact local step, and shows that both satisfy \eqref{eq:comparison}.

\begin{theorem}[What the Frank--Wolfe step keeps]\label{thm:fw-approx}
Let Assumption~\ref{ass:set} hold, fix $x\in\X$, $g\in\R^n$, $\gamma\in[0,1]$, and $t=\gamma D$, and let $y=x+\gamma(s-x)$ with $s$ from \eqref{eq:oracle}. Then $y$ satisfies the contracted comparison \eqref{eq:comparison}, and
\begin{equation}\label{eq:model-fraction}
 \gamma M_t(x)\le\ip g{x-y}=\gamma G\le M_t(x).
\end{equation}
The exact local step $b$ of \eqref{eq:brox} also satisfies \eqref{eq:comparison}, with $b$ in place of $y$, and $\ip g{x-b}=M_t(x)$.
\end{theorem}
\begin{proof}
Since $\mathcal C_t(x)\subseteq\X$,
\[
 M_t(x)=\max_{z\in\mathcal C_t(x)}\ip g{x-z}
 \le\max_{z\in\X}\ip g{x-z}=G,
\]
hence $\gamma M_t(x)\le\gamma G=\ip g{x-y}$.
On the other hand, by \eqref{eq:contracted-value} and \eqref{eq:inclusion},
\[
 \gamma G=\ip g{x-y}
 =\max_{z\in\mathcal H_\gamma(x)}\ip g{x-z}
 \le\max_{z\in\mathcal C_t(x)}\ip g{x-z}=M_t(x).
\]
Both $y$ and $b$ achieve a decrease of at least $\gamma G$.
Consequently, for either point $v\in\{y,b\}$,
\[
 \ip g{v-x}\le-\gamma G
 =\gamma\min_{z\in\X}\ip g{z-x}
 \le\gamma\ip g{z-x}\qquad(z\in\X),
\]
which is \eqref{eq:comparison}.
\end{proof}

When $M_t(x)>0$, the ratio $\ip g{x-y}/M_t(x)$ is the fraction of the exact local decrease attained by Frank--Wolfe. By \eqref{eq:model-fraction}, this fraction is at least $\gamma$. The local-model error satisfies
\begin{equation}\label{eq:fw-local-error}
 0\le M_t(x)-\ip g{x-y}\le(1-\gamma)M_t(x).
\end{equation}
The factor $\gamma$ cannot be improved uniformly. On $\X=[0,D]$ with $g=-1$, $0<\gamma\le1$, and $x=(1-\gamma)D$, the global minimizer $s=D$ lies within the radius $t=\gamma D$ of $x$, so $M_t(x)=\gamma D$, while the Frank--Wolfe step attains $\gamma^2D$: exactly the fraction $\gamma$.

The contraction is the smallest compact convex region whose linear minimizer satisfies the contracted comparison for every linear objective. Here the region must be chosen independently of $g$.

\begin{proposition}[The smallest region preserving the comparison]\label{prop:minimality}
Let Assumption~\ref{ass:set} hold, and fix $x\in\X$ and $\gamma\in[0,1]$. Let $\mathcal R\subseteq\mathcal C_{\gamma D}(x)$ be nonempty, compact, and convex, chosen independently of $g$. Then
\begin{equation}\label{eq:uniform-comparison}
 \min_{u\in\mathcal R}\ip g{u-x}\le\gamma\min_{z\in\X}\ip g{z-x}\qquad\text{for every }g\in\R^n
\end{equation}
holds if and only if $\mathcal H_\gamma(x)\subseteq\mathcal R$. Thus $\mathcal H_\gamma(x)$ is the smallest such region by inclusion.
\end{proposition}
\begin{proof}
If $\mathcal H_\gamma(x)\subseteq\mathcal R$, the minimum over $\mathcal R$ is at most the minimum over $\mathcal H_\gamma(x)$, which equals the right side of \eqref{eq:uniform-comparison} by \eqref{eq:contracted-value}. Conversely, suppose some $\tilde z\in\mathcal H_\gamma(x)$ lies outside $\mathcal R$. Let $\bar z$ be the point of $\mathcal R$ closest to $\tilde z$; it exists and is unique because $\mathcal R$ is compact and convex, and $\bar z\ne\tilde z$.
By the first-order optimality condition for the projection of $\tilde z$ onto $\mathcal R$,
\[
 \ip{\bar z-\tilde z}{u-\bar z}\ge0\qquad(u\in\mathcal R).
\]
With $g\eqdef\bar z-\tilde z$, every $u\in\mathcal R$ satisfies
\[
 \ip g{u-\tilde z}=\ip{\bar z-\tilde z}{u-\bar z}+\norm{\bar z-\tilde z}^2\ge\norm{\bar z-\tilde z}^2>0.
\]
Hence
\[
 \min_{u\in\mathcal R}\ip g{u-x}>\ip g{\tilde z-x}\ge\min_{v\in\mathcal H_\gamma(x)}\ip g{v-x}
 \overset{\eqref{eq:contracted-value}}{=}\gamma\min_{z\in\X}\ip g{z-x},
\]
which contradicts \eqref{eq:uniform-comparison} for this $g$.
\end{proof}
The uniformity in $g$ is essential: the proposition characterizes convex regions chosen before the linear objective is known. It does not assert minimality among direction-dependent constructions.

\paragraph{The fraction alone does not imply convergence.}
The sandwich \eqref{eq:model-fraction} raises the question whether the weaker bound $\ip g{x-y}\ge\gamma M_t(x)$, rather than the contracted comparison $\ip g{x-y}\ge\gamma G$, suffices for convergence. The following sequence satisfies the weaker bound with equality at every step but converges to a nonoptimal point. On $\X=[0,1]$ with $f(x)=x$, so that $g=1$ and $D=1$, let $\gamma_k=t_k=1/(k+2)$, $\widehat x_0=1$, and $\widehat x_{k+1}=\widehat x_k-\gamma_k^2$. After $k$ updates, with step fractions $1/2,\dots,1/(k+1)$,
\[
 \widehat x_k=1-\sum_{j=2}^{k+1}\frac1{j^2}\ \ge\ 1-\sum_{j=2}^{k+1}\Bigl(\frac1{j-1}-\frac1j\Bigr)=\frac1{k+1}>t_k.
\]
Since $\widehat x_k-t_k\in[0,1]$, the linear objective is minimized over the ball intersection at $\widehat x_k-t_k$, so $M_{t_k}(\widehat x_k)=t_k$. Thus the step attains $\widehat x_k-\widehat x_{k+1}=\gamma_k^2=\gamma_kM_{t_k}(\widehat x_k)$. Yet $\widehat x_k\to2-\pi^2/6>0$. The sequence is not a Frank--Wolfe sequence: the oracle returns $s_k=0$, and \eqref{eq:comparison} at $z=0$ requires $\widehat x_{k+1}-\widehat x_k\le-\gamma_k\widehat x_k$, a decrease of at least $\gamma_k\widehat x_k>\gamma_k^2$. Frank--Wolfe with the same step fractions gives $x_{k+1}=(1-\gamma_k)x_k$ and $x_k=1/(k+1)\to0$.

\paragraph{What the approximation gains on the simplex.}
The restriction that costs local accuracy can preserve sparsity. The example compares the number of nonzero coordinates of the exact local step $b$ and of a Frank--Wolfe step $y$ at the same radius. Let $\X=\{z\in\R^m_+:\sum_iz_i=1\}$ be the probability simplex, $m\ge3$, so that $D=\sqrt2$, and take $x=e_1$, $g=e_1$, and $0<t<1$. Since $\ip g z=z_1$, write $z_1=1-\beta$, so that $\sum_{i=2}^m z_i=\beta$.
Minimizing $z_1$ amounts to maximizing $\beta$, subject to the ball constraint $\norm{z-e_1}\le t$.
By Cauchy--Schwarz,
\begin{equation}\label{eq:simplex-radius}
 t^2\ge\norm{z-e_1}^2
 =\beta^2+\sum_{i=2}^m z_i^2
 \ge\beta^2+\frac{\left(\sum_{i=2}^m z_i\right)^2}{m-1}
 =\frac m{m-1}\beta^2.
\end{equation}
Hence $\beta\le\beta_t\eqdef t\sqrt{(m-1)/m}$.
This bound is attained when $z_2=\cdots=z_m=\beta_t/(m-1)$, giving the exact local step.
A vertex-returning global oracle selects $e_j$ with $j\ge2$, say $e_2$, and with $\gamma=t/\sqrt2$ the two steps are
\begin{equation}\label{eq:simplex-steps}
 b=(1-\beta_t)e_1+\frac{\beta_t}{m-1}\sum_{i=2}^m e_i,\qquad
 y=(1-\gamma)e_1+\gamma e_2.
\end{equation}
Both move the full distance $t$ from $x$, but $b$ has $m$ nonzero coordinates and $y$ has two. The decrease retained is substantial:
\begin{equation}\label{eq:simplex-fraction}
 M_t(e_1)=\beta_t,\qquad
 \frac{\ip g{x-y}}{M_t(x)}=\frac{\gamma}{\beta_t}=\sqrt{\frac{m}{2(m-1)}}\ \ge\ \frac1{\sqrt2}.
\end{equation}

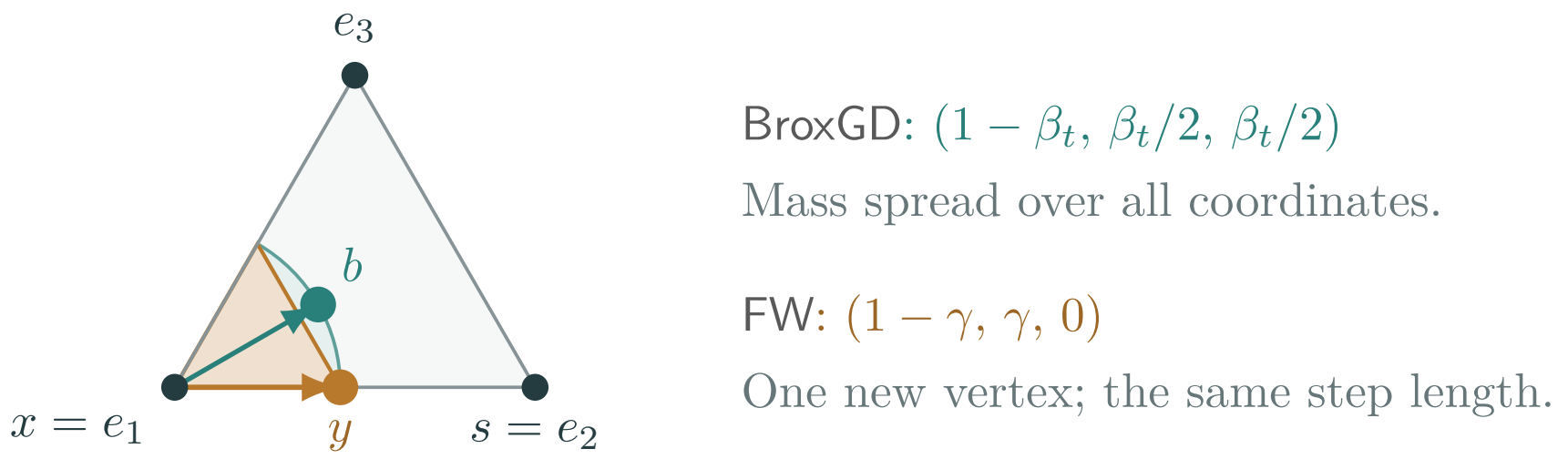
\begin{figure}[!htbp]
\centering
\begin{tikzpicture}[scale=1.85,>=Latex,font=\small]
% An isometric planar embedding of the three-dimensional probability simplex.
% Edges have length sqrt(2); t=.65, b=(t cos(30),t sin(30)), y=(t,0).
\coordinate (sx) at (0,0);
\coordinate (s2) at ({sqrt(2)},0);
\coordinate (s3) at ({sqrt(2)/2},{sqrt(6)/2});
\coordinate (sb) at ({.65*sqrt(3)/2},.325);
\coordinate (sy) at (.65,0);
\fill[ink!4] (sx)--(s2)--(s3)--cycle;
\fill[jade!12] (sx)--(.65,0) arc[start angle=0,end angle=60,radius=.65]--cycle;
\draw[jade!75,line width=.7pt] (.65,0) arc[start angle=0,end angle=60,radius=.65];
\filldraw[fill=ochre!23,draw=ochre,line width=.8pt]
 (sx)--(sy)--(.325,{.65*sqrt(3)/2})--cycle;
\draw[ink!55,line width=.7pt] (sx)--(s2)--(s3)--cycle;
\draw[->,jade,line width=1pt] (sx)--(sb);
\draw[->,ochre,line width=1.1pt] (sx)--(sy);
\fill[ink] (sx) circle (1.5pt) node[below left=2pt] {$x=e_1$};
\fill[ink] (s2) circle (1.5pt) node[below=3pt] {$s=e_2$};
\fill[ink] (s3) circle (1.5pt) node[above=3pt] {$e_3$};
\fill[jade] (sb) circle (2pt) node[above right=1pt,text=jade] {$b$};
\fill[ochre] (sy) circle (2pt) node[below=3pt,text=ochre!85!black] {$y$};
\node[anchor=west,text=jade] at (2.15,1.02) {\algname{BroxGD}: $(1-\beta_t,\,\beta_t/2,\,\beta_t/2)$};
\node[anchor=west,text=ink!70] at (2.15,.72) {Mass spread over all coordinates.};
\node[anchor=west,text=ochre!85!black] at (2.15,.27) {\algname{FW}: $(1-\gamma,\,\gamma,\,0)$};
\node[anchor=west,text=ink!70] at (2.15,-.03) {One new vertex; the same step length.};
\end{tikzpicture}
\caption{On the simplex with $m=3$, the \algname{BroxGD} step $b$ has three nonzero coordinates and the \algname{FW} step $y$ has two. Both move the same Euclidean distance $t$.}
\label{fig:simplex}
\end{figure}

Figure~\ref{fig:simplex} illustrates the case $m=3$. More generally, a vertex oracle adds at most one new vertex per step, the sparsity property of Frank--Wolfe \citep{jaggi}.

\section{One convergence proof for both methods}\label{sec:convergence}
The convergence proof of this section uses only two inequalities per step: a bound on the step length, and the comparison \eqref{eq:comparison} up to an additive error term. Every method whose iterates meet both, exact \algname{BroxGD} and Frank--Wolfe included, satisfies the same recurrence for the objective gap.

\begin{assumption}[Smooth convex objective]\label{ass:objective}
The function $f$ is differentiable on an open neighborhood of $\X$. It is convex on $\X$, that is,
\begin{equation}\label{eq:convex}
 f(z)\ge f(x)+\ip{\nabla f(x)}{z-x}\qquad(x,z\in\X),
\end{equation}
and $L$-smooth on $\X$ for some $L>0$, that is,
\begin{equation}\label{eq:smooth}
 f(z)\le f(x)+\ip{\nabla f(x)}{z-x}+\frac L2\norm{z-x}^2\qquad(x,z\in\X).
\end{equation}
\end{assumption}
Under Assumptions~\ref{ass:set} and~\ref{ass:objective}, fix a minimizer $x_\star$ and set
\begin{equation}\label{eq:notation}
 f_\star\eqdef f(x_\star),\qquad \Delta(x)\eqdef f(x)-f_\star,\qquad g\eqdef\nabla f(x).
\end{equation}

\begin{theorem}[The comparison that preserves convergence]\label{thm:approx}
Let Assumptions~\ref{ass:set} and~\ref{ass:objective} hold. Let $x_k\in\X$, $\gamma_k\in[0,1]$, and $e_k\ge0$ satisfy, with $g_k\eqdef\nabla f(x_k)$,
\begin{equation}\label{eq:approx-condition}
 \norm{x_{k+1}-x_k}\le\gamma_kD,\qquad
 \ip{g_k}{x_{k+1}-x_k}\le\gamma_k\ip{g_k}{z-x_k}+e_k\quad(z\in\X).
\end{equation}
Then $\Delta_k\eqdef\Delta(x_k)$ satisfies
\begin{equation}\label{eq:recurrence}
 \Delta_{k+1}\le(1-\gamma_k)\Delta_k+\frac{LD^2}{2}\gamma_k^2+e_k.
\end{equation}
If $e_k\le A\gamma_k^2$ for some $A\ge0$ and $\gamma_k=\nicefrac{a}{(k+a)}$ for some $a>1$, then
\begin{equation}\label{eq:rate}
 \Delta_k\le\frac{a^2(LD^2+2A)}{2(a-1)(k+a)}\qquad(k\ge1).
\end{equation}
\end{theorem}
The error $e_k$ measures how far the step falls short of the comparison \eqref{eq:comparison}. It is zero for exact \algname{BroxGD} at radius $\gamma_kD$ and for Frank--Wolfe.
\begin{proof}
Smoothness \eqref{eq:smooth} at $z=x_{k+1}$, after subtracting $f_\star$, gives
\[
 \Delta_{k+1}\overset{\eqref{eq:smooth}}{\le}\Delta_k+\ip{g_k}{x_{k+1}-x_k}+\frac L2\norm{x_{k+1}-x_k}^2.
\]
The comparison at $z=x_\star$ and convexity bound the linear term:
\begin{equation}\label{eq:gap-comparison}
 \ip{g_k}{x_{k+1}-x_k}
 \overset{\eqref{eq:approx-condition}}{\le}\gamma_k\ip{g_k}{x_\star-x_k}+e_k
 \overset{\eqref{eq:convex}}{\le}-\gamma_k\Delta_k+e_k.
\end{equation}
The step bound in \eqref{eq:approx-condition} now gives \eqref{eq:recurrence}.

For the rate, put $C\eqdef LD^2+2A$. The error bound turns \eqref{eq:recurrence} into
\begin{equation}\label{eq:C-recurrence}
 \Delta_{k+1}\le(1-\gamma_k)\Delta_k+\frac C2\gamma_k^2.
\end{equation}
Set $\Lambda\eqdef\nicefrac{Ca^2}{2(a-1)}$. We prove $\Delta_k\le\Lambda/(k+a)$ for $k\ge1$ by induction. Since $\gamma_0=1$, \eqref{eq:C-recurrence} gives $\Delta_1\le C/2\le\Lambda/(a+1)$. If $\Delta_k\le\Lambda/(k+a)$ for some $k\ge1$, then
\[
 \Delta_{k+1}\overset{\eqref{eq:C-recurrence}}{\le}
 \left(1-\frac a{k+a}\right)\frac\Lambda{k+a}+\frac{Ca^2}{2(k+a)^2}
 =\frac\Lambda{k+a}-\frac\Lambda{(k+a)^2}\le\frac\Lambda{k+a+1},
\]
where the equality uses $a\Lambda-Ca^2/2=\Lambda$ and the last inequality $(k+a-1)(k+a+1)\le(k+a)^2$.
\end{proof}

\begin{corollary}[Exact \algname{BroxGD} and Frank--Wolfe share the proof]\label{cor:inherited}
Let Assumptions~\ref{ass:set} and~\ref{ass:objective} hold. Exact \algname{BroxGD} at radius $t_k=\gamma_kD$ and Frank--Wolfe with step fraction $\gamma_k$ both satisfy
\begin{equation}\label{eq:exact-recurrence}
 \Delta_{k+1}\le(1-\gamma_k)\Delta_k+\frac{LD^2}{2}\gamma_k^2.
\end{equation}
With $\gamma_k=\nicefrac{2}{(k+2)}$, both satisfy
\begin{equation}\label{eq:exact-rate}
 \Delta_k\le\frac{2LD^2}{k+2}\qquad(k\ge1).
\end{equation}
\end{corollary}
\begin{proof}
Both steps have length at most $\gamma_kD$, by the radius for \algname{BroxGD} and by Lemma~\ref{lem:contraction} for Frank--Wolfe, and both satisfy \eqref{eq:comparison} by Theorem~\ref{thm:fw-approx}. So \eqref{eq:approx-condition} holds with $e_k=0$, and Theorem~\ref{thm:approx} with $A=0$ and $a=2$ gives the claims.
\end{proof}

For \algname{BroxGD}, the diameter dependence in Corollary~\ref{cor:inherited} comes from the radius schedule $t_k=\gamma_kD$. Other radius choices give diameter-independent guarantees, as shown by the fixed-horizon and backtracking results of \citet[Theorems~3.2 and~3.6]{broxgd}. These results also apply on unbounded sets; the former uses an initial-distance bound, the latter a distance bound along the accepted trajectory.

Every $a>1$ gives the same order $O(1/k)$; the constant $a^2/(a-1)=4+(a-2)^2/(a-1)$ is smallest at the familiar choice $a=2$.

The same argument covers exact minimization over any compact region between $\mathcal H_{\gamma_k}(x_k)$ and $\mathcal C_{\gamma_kD}(x_k)$, and inexact solves over such regions. If the global oracle returns $s_k\in\X$ with
\begin{equation}\label{eq:global-error}
 \ip{g_k}{s_k}\le\min_{z\in\X}\ip{g_k}{z}+\varepsilon_k,
\end{equation}
the Frank--Wolfe step satisfies \eqref{eq:approx-condition} with $e_k=\gamma_k\varepsilon_k$, so $\varepsilon_k\le A\gamma_k$ suffices for \eqref{eq:rate}; this is the approximate-oracle analysis of \citet{jaggi}. Each Frank--Wolfe step uses one gradient and one global oracle call, and its prescribed schedule does not require $D$.

The schedule $\gamma_k=a/(k+a)$ gives the stated convergence rate. For convergence alone, the step fractions can satisfy more general conditions.

\begin{corollary}[Diminishing, nonsummable step fractions]\label{cor:diminishing}
Let Assumptions~\ref{ass:set} and~\ref{ass:objective} hold, and let the iterates satisfy \eqref{eq:approx-condition} with $e_k\le A\gamma_k^2$ for a fixed $A\ge0$. If $\gamma_k\in[0,1]$, $\gamma_k\to0$, and $\sum_{k=0}^{\infty}\gamma_k=\infty$, then $\Delta_k\to0$. In particular this holds for exact \algname{BroxGD} at radius $\gamma_kD$ and for Frank--Wolfe with step fraction $\gamma_k$.
\end{corollary}
\begin{proof}
Put $\kappa\eqdef LD^2/2+A$, so that \eqref{eq:recurrence} reads $\Delta_{k+1}\le(1-\gamma_k)\Delta_k+\kappa\gamma_k^2$.
Fix $\epsilon>0$. Since $\gamma_k\to0$, choose $N$ with $\kappa\gamma_k\le\epsilon$ for all $k\ge N$.
Then $\kappa\gamma_k^2\le\epsilon\gamma_k$, and
\[
 \Delta_{k+1}-\epsilon
 \le(1-\gamma_k)\Delta_k+\epsilon\gamma_k-\epsilon
 =(1-\gamma_k)(\Delta_k-\epsilon)\qquad(k\ge N).
\]
Since $1-\gamma_k\ge0$, iterating gives
\[
 \Delta_k-\epsilon
 \le(\Delta_N-\epsilon)\prod_{j=N}^{k-1}(1-\gamma_j).
\]
The product tends to zero because
\[
 0\le\prod_{j=N}^{k-1}(1-\gamma_j)
 \le\exp\Bigl(-\sum_{j=N}^{k-1}\gamma_j\Bigr)\to0,
\]
where the last limit follows from $\sum_k\gamma_k=\infty$.
Hence $\limsup_k\Delta_k\le\epsilon$.
Since $\epsilon>0$ is arbitrary and $\Delta_k\ge0$, we obtain $\Delta_k\to0$.
\end{proof}
For Frank--Wolfe this is the open-loop result of \citet{dunn}. No square-summability is needed: $\gamma_k=(k+1)^{-p}$ works for every $0<p\le1$, including $p\le\nicefrac12$, where $\sum_k\gamma_k^2$ diverges. At $p=1$ the recurrence gives the more specific bound
\begin{equation}\label{eq:harmonic-rate}
 \Delta_k\le\frac{\kappa H_k}{k},\qquad
 H_k\eqdef\sum_{j=1}^{k}\frac1j\qquad(k\ge1),
\end{equation}
since multiplying \eqref{eq:recurrence} by $k+1$ gives $(k+1)\Delta_{k+1}\le k\Delta_k+\kappa/(k+1)$, and summing telescopes. This gives an $O(\log k/k)$ bound for $\gamma_k=\nicefrac1{k+1}$.

We now return to the two approximations used to obtain Frank--Wolfe from \algname{BPM} and compare their effect on the objective value after one step. Fix $x$ and $t=\gamma D$, and let $p$ minimize $f$ over $\mathcal C_t(x)$, giving the \algname{BPM} step.

\begin{lemma}[The two approximation errors]\label{lem:two-errors}
Let Assumptions~\ref{ass:set} and~\ref{ass:objective} hold. Fix $x\in\X$, $\gamma\in(0,1]$, $t=\gamma D$, and $g=\nabla f(x)$. Let $p\in\argmin_{z\in\mathcal C_t(x)}f(z)$ and let $y$ be the Frank--Wolfe step \eqref{eq:fw-step}. Then
\begin{equation}\label{eq:two-errors}
 0\le f(y)-f(p)\le
 \underbrace{M_t(x)-\gamma G}_{\text{local-model error}}
 +\underbrace{\frac L2\norm{y-x}^2}_{\text{linearization bound}}.
\end{equation}
For the exact \algname{BroxGD} step $b$, the local-model error vanishes and $0\le f(b)-f(p)\le\nicefrac{Lt^2}{2}$.
\end{lemma}
\begin{proof}
For any $v\in\mathcal C_t(x)$, applying \eqref{eq:smooth} with $z=v$ and \eqref{eq:convex} with $z=p$, both based at $x$, gives
\begin{equation}\label{eq:ball-objective-error}
 f(v)-f(p)\overset{\eqref{eq:smooth},\eqref{eq:convex}}{\le}\ip g{v-p}+\frac L2\norm{v-x}^2
 \overset{\eqref{eq:local-value}}{\le} M_t(x)-\ip g{x-v}+\frac L2\norm{v-x}^2,
\end{equation}
where the second step writes $\ip g{v-p}=\ip g{x-p}-\ip g{x-v}$ and uses $\ip g{x-p}\le M_t(x)$, since $p\in\mathcal C_t(x)$. Taking $v=y$ and $\ip g{x-y}=\gamma G$ from Theorem~\ref{thm:fw-approx} gives \eqref{eq:two-errors}. Taking $v=b$, $\ip g{x-b}=M_t(x)$ and $\norm{b-x}\le t$. The lower bounds hold because $p$ minimizes $f$ over $\mathcal C_t(x)$, which contains $y$ and $b$.
\end{proof}

\algname{BroxGD} therefore approximates \algname{BPM} in objective value with an error bounded by $Lt^2/2$. The loss from restricting the region in Frank--Wolfe can instead be proportional to $t$, as the triangle example in Section~\ref{sec:geometry} shows.

\section{Choosing the step size along the oracle direction}\label{sec:radius}
Section~\ref{sec:convergence} used prescribed step fractions. We now choose the step after the oracle call by minimizing the smoothness bound along the segment from $x$ to $s$. This gives the standard short-step Frank--Wolfe rule, which retains the $O(1/k)$ convergence bound and guarantees an objective decrease at every nonoptimal iterate. Prescribed step fractions can instead increase the objective, as the triangle example in Section~\ref{sec:geometry} shows.

With $s$ and $G$ from \eqref{eq:oracle}, convexity gives
\begin{equation}\label{eq:gap}
 G\overset{\eqref{eq:oracle}}{\ge}\ip g{x-x_\star}
 \overset{\eqref{eq:convex}}{\ge}\Delta(x),
\end{equation}
so $G=0$ means that $x$ is optimal. Otherwise $d\eqdef s-x\ne0$, and we write
\begin{equation}\label{eq:direction}
 T\eqdef\norm d,\qquad u\eqdef d/T,\qquad q\eqdef G/T=-\ip g u.
\end{equation}
Here $T$ is the distance to the oracle point, $u$ is the unit direction toward it, and $q$ is the linear-model decrease per unit distance.

For $0\le\tau\le T$, convexity ensures that $x+\tau u$ is feasible.
This step has length $\tau$ and decreases the linear model by $\tau q$.
Smoothness therefore gives
\begin{equation}\label{eq:segment-bound}
 f(x+\tau u)\overset{\eqref{eq:smooth},\eqref{eq:direction}}{\le}
 f(x)-\tau q+\frac L2\tau^2\qquad(0\le\tau\le T).
\end{equation}
The Frank--Wolfe step with fraction $\gamma$ corresponds to $\tau=\gamma T$.
This actual step length can be smaller than the comparison radius $t=\gamma D$, which bounds the entire contracted set.

We choose the step length by minimizing the quadratic upper bound in \eqref{eq:segment-bound} over $0\le\tau\le T$.

\begin{theorem}[Short-step Frank--Wolfe from step-length selection]\label{thm:short}
Let Assumptions~\ref{ass:set} and~\ref{ass:objective} hold, let $G>0$, and choose
\begin{equation}\label{eq:radius}
 \tau^+\eqdef\min\{T,\,q/L\},\qquad x^+\eqdef x+\tau^+u.
\end{equation}
This is the short-step Frank--Wolfe update, since
\begin{equation}\label{eq:short-step}
 \gamma\eqdef\frac{\tau^+}{T}
 \overset{\eqref{eq:radius},\eqref{eq:direction}}{=}\min\left\{1,\frac{G}{L\norm d^2}\right\},
 \qquad x^+=x+\gamma d.
\end{equation}
It satisfies
\begin{equation}\label{eq:descent}
 f(x)-f(x^+)\ge\frac12\min\{G,\,q^2/L\}.
\end{equation}
Starting from any $x_0\in\X$ and stopping when $G=0$, the iterates satisfy $\Delta_k\le2LD^2/(k+2)$ for $k\ge1$, with the sequence held constant after termination.
\end{theorem}
\begin{proof}
The right side of \eqref{eq:segment-bound} is a convex quadratic in $\tau$ with unconstrained minimizer $q/L$, so \eqref{eq:radius} minimizes it over $[0,T]$.

If $q/L\le T$, then $\tau^+=q/L$, and \eqref{eq:segment-bound} gives
\[
 f(x)-f(x^+)\ge\frac{q^2}{2L}.
\]
Otherwise $q/L>T$, so $LT<q$. At $\tau^+=T$, the decrease satisfies
\[
 f(x)-f(x^+)\ge qT-\frac L2T^2
 \ge\frac{qT}{2}=\frac G2.
\]
Together the two cases give \eqref{eq:descent}.

For the rate, we compare the chosen step fraction $\gamma$ with an arbitrary fraction $\eta\in[0,1]$.
Write the smoothness bound as
\[
 \varphi(\eta)\eqdef\Delta(x)-\eta G+\frac L2\eta^2\norm d^2\qquad(0\le\eta\le1),
\]
which is \eqref{eq:segment-bound} with $\tau=\eta T$ after subtracting $f_\star$.
The rule \eqref{eq:short-step} minimizes $\varphi$ over $[0,1]$, so
\[
 \Delta(x^+)\le\varphi(\gamma)\le\varphi(\eta)\qquad(0\le\eta\le1).
\]
Using $G\ge\Delta(x)$ from \eqref{eq:gap} and $\norm d\le D$,
\[
 \Delta(x^+)\le(1-\eta)\Delta(x)+\frac{LD^2}{2}\eta^2\qquad(0\le\eta\le1).
\]
Taking $\eta=\nicefrac{2}{(k+2)}$ at iteration $k$ gives the recurrence \eqref{eq:C-recurrence} with $C=LD^2$ and $a=2$, and the induction in the proof of Theorem~\ref{thm:approx} gives the bound.
\end{proof}

Exact line search minimizes $f$ over the same segment, so from any current point it decreases the objective at least as much as the short step. It therefore satisfies the sublinear bound of Theorem~\ref{thm:short} and the interior linear bound of Corollary~\ref{cor:interior} below, at the cost of a one-dimensional minimization.

We can obtain a linear rate by showing that each step decreases the objective by at least a fixed positive fraction of the current objective gap. The following additional assumptions allow us to establish such a bound.

\begin{assumption}[Strong convexity and an interior minimizer]\label{ass:interior}
The function $f$ is $\mu$-strongly convex on $\X$, with $\mu>0$, and some minimizer $x_\star$ satisfies $\B(x_\star,r)\subseteq\X$ for some $r>0$.
\end{assumption}

The interior ball bounds the oracle slope from below in terms of the gradient norm, while strong convexity relates the gradient norm to the objective gap. Together, these bounds strengthen \eqref{eq:descent} to a linear rate.

\begin{corollary}[Interior-solution linear convergence]\label{cor:interior}
Let Assumptions~\ref{ass:set} to~\ref{ass:interior} hold. Short-step Frank--Wolfe satisfies
\begin{equation}\label{eq:interior-rate}
 \Delta_k\le(1-\rho)^k\Delta_0,\qquad
 \rho\eqdef\min\left\{\frac12,\frac{\mu r^2}{LD^2}\right\}.
\end{equation}
This is the classical interior-solution regime of \citet{guelat}.
\end{corollary}
\begin{proof}
At a nonoptimal $x$, convexity gives $g\ne0$. The point $z=x_\star-rg/\norm g$ is feasible by Assumption~\ref{ass:interior}, so
\begin{equation}\label{eq:interior-gap}
 G\overset{\eqref{eq:oracle}}{\ge}\ip g{x-z}
 =\ip g{x-x_\star}+r\norm g
 \overset{\eqref{eq:convex}}{\ge}\Delta(x)+r\norm g\ge r\norm g.
\end{equation}
Strong convexity at $x$, evaluated at $x_\star$, gives with $\delta\eqdef\norm{x-x_\star}$
\begin{equation}\label{eq:gradient-gap}
 \Delta(x)\le\ip g{x-x_\star}-\frac\mu2\delta^2
 \le\norm g\,\delta-\frac\mu2\delta^2
 =\frac{\norm g^2}{2\mu}-\frac\mu2\Bigl(\delta-\frac{\norm g}{\mu}\Bigr)^2
 \le\frac{\norm g^2}{2\mu}.
\end{equation}
Hence
\begin{equation}\label{eq:interior-slope}
 q^2\overset{\eqref{eq:direction}}{=}\frac{G^2}{T^2}
 \overset{\eqref{eq:diameter}}{\ge}\frac{G^2}{D^2}
 \overset{\eqref{eq:interior-gap}}{\ge}\frac{r^2\norm g^2}{D^2}
 \overset{\eqref{eq:gradient-gap}}{\ge}\frac{2\mu r^2}{D^2}\Delta(x).
\end{equation}
Inserting \eqref{eq:gap} and \eqref{eq:interior-slope} into \eqref{eq:descent},
\[
 f(x)-f(x^+)\ge\frac12\min\left\{\Delta(x),\frac{2\mu r^2}{LD^2}\Delta(x)\right\}=\rho\,\Delta(x),
\]
so $\Delta(x^+)\le(1-\rho)\Delta(x)$, and iteration gives \eqref{eq:interior-rate}. If an optimum is reached, the sequence is held constant.
\end{proof}

The ball assumption requires $\X$ to be full-dimensional. The same rate holds for a lower-dimensional set if the ball condition is satisfied within its affine hull. Let $\mathcal T$ be the direction space of that hull and $g_{\mathcal T}$ the projection of $g$ onto $\mathcal T$. Feasible displacements lie in $\mathcal T$, so this replacement leaves the Frank--Wolfe gap unchanged. At a nonoptimal point $g_{\mathcal T}\ne0$, and the proof gives $G\ge r\norm{g_{\mathcal T}}$ and $\Delta(x)\le\norm{g_{\mathcal T}}^2/(2\mu)$, yielding \eqref{eq:interior-rate} with the same $\rho$.

\section{What a coarse approximation can lose}\label{sec:geometry}
The previous sections showed that Frank--Wolfe preserves the comparison used in the objective convergence proof. We now show that, at a radius where \algname{BroxGD} guarantees a decrease in distance to the minimizer, Frank--Wolfe can move farther away.

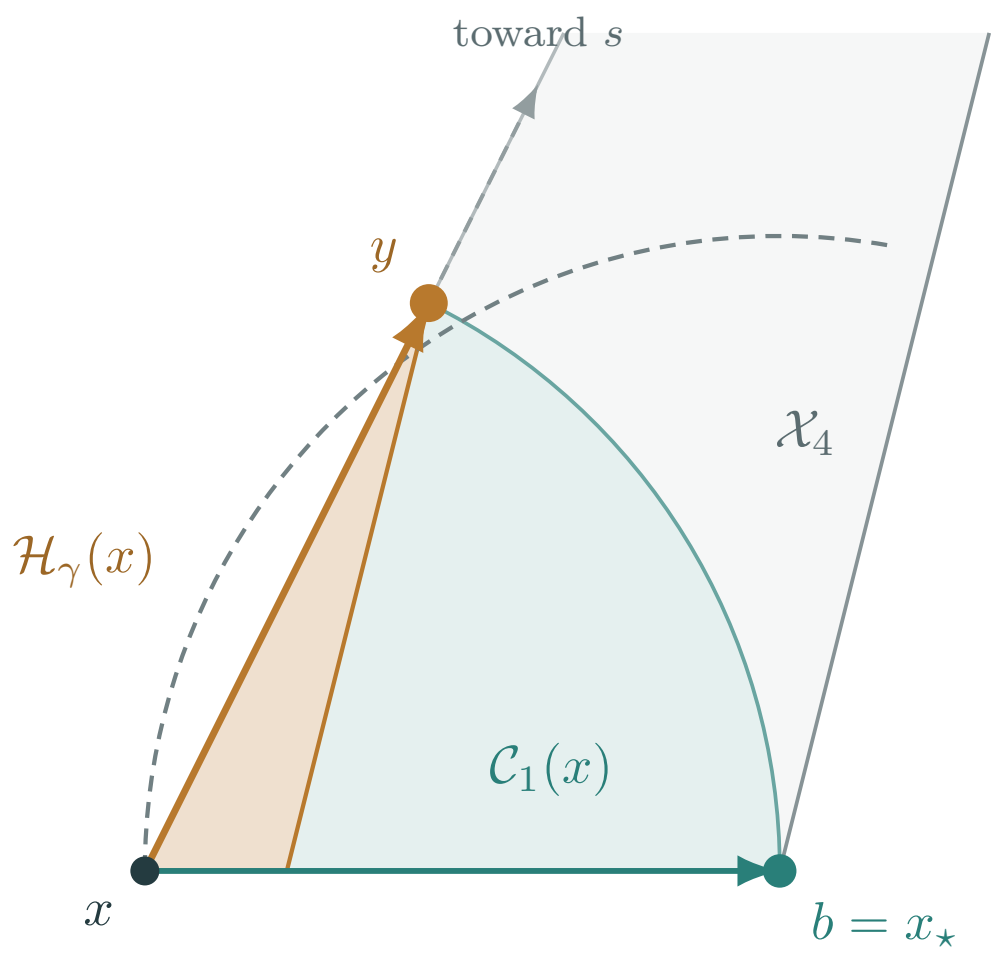
\begin{figure}[!htbp]
\centering
\begin{tikzpicture}[scale=4.2,>=Latex,font=\small]
% Close-up of X_4 near the unit ball, with equal Euclidean axis scales.
% The distant oracle vertex (2,4) is indicated by the continuing direction.
\coordinate (x) at (0,0);
\coordinate (xs) at (1,0);
\coordinate (y) at ({1/sqrt(5)},{2/sqrt(5)});
\coordinate (h) at ({1/sqrt(20)},0);
\coordinate (ul) at (.66,1.32);
\coordinate (ur) at (1.33,1.32);
\pgfmathsetmacro{\rayangle}{atan(2)}
% The top edge is omitted because the feasible triangle continues upward.
\fill[ink!4] (x)--(xs)--(ur)--(ul)--cycle;
\fill[jade!12] (x)--(1,0)
 arc[start angle=0,end angle=\rayangle,radius=1]--cycle;
\draw[ink!55,line width=.7pt] (x)--(xs)--(ur);
\draw[ink!35,line width=.65pt] (x)--(ul);
\draw[jade!70,line width=.7pt] (1,0)
 arc[start angle=0,end angle=\rayangle,radius=1];
\filldraw[fill=ochre!23,draw=ochre,line width=.8pt]
 (x)--(h)--(y)--cycle;
% The comparison circle has center x_star and radius ||x-x_star||=1.
\draw[ink!65,densely dashed,line width=.8pt]
 (x) arc[start angle=180,end angle=80,radius=1];
\draw[->,ink!50,densely dashed,line width=.75pt]
 (y)--(.62,1.24)
 node[above=3pt,font=\scriptsize,text=ink!75] {toward $s$};
\draw[->,jade,line width=1.1pt] (x)--(xs);
\draw[->,ochre,line width=1.3pt] (x)--(y);
\fill[ink] (x) circle (.65pt) node[below left=2pt] {$x$};
\fill[jade] (xs) circle (.75pt)
 node[below right=2pt,text=jade] {$b=x_\star$};
\fill[ochre] (y) circle (.85pt)
 node[above left=2pt,text=ochre!85!black] {$y$};
\node[text=jade] at (.64,.16) {$\mathcal C_1(x)$};
\node[anchor=east,text=ochre!85!black] at (.05,.49)
 {$\mathcal H_\gamma(x)$};
\node[text=ink!75] at (1.04,.69) {$\X_4$};
\end{tikzpicture}
\caption{Close-up of the triangle example for $h=4$ and $t=1$. The dashed circle is centered at $x_\star$ and passes through $x$; the \algname{FW} step $y$ lies outside it, whereas the exact \algname{BroxGD} step reaches $b=x_\star$.}
\label{fig:triangle}
\end{figure}

Let $h>0$ and
\begin{equation}\label{eq:triangle}
 \X_h\eqdef\operatorname{conv}\{(0,0),(1,0),(2,h)\},\qquad
 f(z)\eqdef\frac12\norm{z-(1,0)}^2.
\end{equation}
Take the current point $x=(0,0)$. The minimizer is $x_\star=(1,0)$, the gradient is $g=(-1,0)$, and the diameter is $D_h=\sqrt{4+h^2}$.
Figure~\ref{fig:triangle} shows the case $h=4$.

Take $t=1$. Since $\ip g{z-x}=-z_1$, the local problem maximizes the first coordinate.
Over the unit ball, $z_1\le1$, with equality only at $(1,0)$.
This point belongs to $\X_h$, so the \algname{BroxGD} step is $b=(1,0)=x_\star$, with linear-model decrease $M_1(x)=1$.

The global oracle returns $s=(2,h)$. With $\gamma=\nicefrac1{D_h}$, the Frank--Wolfe step is $y=(2,h)/D_h$, and
\begin{equation}\label{eq:triangle-model}
 \ip g{x-y}=\frac2{D_h},\qquad\norm{y-x}=1.
\end{equation}
Both steps have length one. As $h$ grows, the Frank--Wolfe step spends most of its length moving vertically, and the fraction $2/D_h$ of the local decrease that it attains tends to zero.

The same example illustrates the different dependence on the radius of the two errors in Lemma~\ref{lem:two-errors}.
For $0<t\le1$ and $\gamma=\nicefrac t{D_h}$, the steps are $b=(t,0)$ and $y=t(2,h)/D_h$, giving
\[
 M_t(x)-\gamma G=t\left(1-\frac2{D_h}\right).
\]
This local-model error is linear in $t$, whereas the linearization bound $Lt^2/2$ is quadratic.

For a \algname{BroxGD} step $b$ with $\nabla f(x)\ne0$, \citet[Theorem~3.1]{broxgd} prove that the radius condition
\begin{equation}\label{eq:type-one}
 0<t\le\frac{\ip{\nabla f(x)}{x-x_\star}}{\norm{\nabla f(x)}}
\end{equation}
implies
\[
 \norm{b-x_\star}^2\le\norm{x-x_\star}^2-t^2.
\]
They call \eqref{eq:type-one} Type-I admissibility of the radius.
For every $h>0$, the right side of \eqref{eq:type-one} equals one, so the condition allows all radii $0<t\le1$.

For the Frank--Wolfe step at radius $t$, however,
\begin{equation}\label{eq:triangle-distance}
 \norm{y-x_\star}^2=1-\frac{4t}{D_h}+t^2.
\end{equation}
Thus, for any fixed $t\in(0,1]$, choosing $D_h>4/t$ makes Frank--Wolfe increase the distance, while \algname{BroxGD} still guarantees a decrease. The objective also increases, since $f(z)=\frac12\norm{z-x_\star}^2$.

These steps use the prescribed fraction $\gamma=\nicefrac t{D_h}$, rather than the short-step choice of Section~\ref{sec:radius}.

\section{Recovering the distance guarantee}\label{sec:extension}
Section~\ref{sec:geometry} shows that restricting the local problem to the contracted set can lose the distance guarantee of \algname{BroxGD}. To recover it, Appendix~\ref{sec:geometric-appendix} allows a larger region inside the ball and improves the local solution when necessary.

On a polytope, the procedure starts with the segment joining the current point to the oracle vertex (Appendix~\ref{sec:refine}). It minimizes the linear model over the part of this segment inside the ball. An additional global oracle call checks an optimality condition for the full local problem. If the check fails, the call supplies a vertex outside the current hull. That vertex is added, and the linear model is minimized again over the intersection of the enlarged hull with the ball.

Appendix~\ref{sec:geometric-appendix} measures the violation detected by this check, calling it an optimality residual. It gives an accuracy threshold sufficient to recover the distance argument of \algname{BroxGD} and proves that the refinement terminates after finitely many vertex additions (Theorem~\ref{thm:threepoint}, Proposition~\ref{prop:refine}). When $t=\gamma D$ with $0<\gamma\le1$, every intermediate candidate retains the contracted comparison used in Section~\ref{sec:convergence}.

For a prescribed iteration budget $K$, choosing the constant radius $t=R/\sqrt{(1-\theta)K}$, where $R>0$ bounds the initial distance $\norm{x_0-x_\star}$ and $\theta\in[0,1)$ is the threshold parameter, gives a bound of $LR^2/(2(1-\theta)K)$ on the smallest objective gap over the run (Corollary~\ref{cor:certified}). Unlike Corollary~\ref{cor:inherited}, this objective-gap bound after $K$ outer steps does not involve the diameter.

For the distance-increasing triangle examples, the check selects the missing vertex $(1,0)$. Adding it recovers the whole triangle, so the next restricted solve gives the exact \algname{BroxGD} step. The additional accuracy requires extra oracle calls and restricted solves.

\section{Relation to the contracted trust region of Nesterov}\label{sec:perspective}
The contracted set was introduced earlier, with a different motivation, by \citet[Section~3]{nesterov2018}, who proposed, for composite objectives $f+\Psi$ with bounded domain $Q=\operatorname{dom}\Psi$ and contraction parameter $\tau_k\in(0,1]$, the step
\[
 x_{k+1}\in\argmin\bigl\{\ip{\nabla f(x_k)}{y}+\Psi(y):\ y=(1-\tau_k)x_k+\tau_k x,\ x\in Q\bigr\},
\]
described it as a trust-region scheme with a linear model whose trust region is a contraction of the feasible set, and noted that it coincides with Frank--Wolfe when $\Psi$ is an indicator function. Its optimality condition is, in that case, the contracted comparison \eqref{eq:comparison}, and the estimating-sequence analysis there uses it in the role that \eqref{eq:comparison} plays in Theorem~\ref{thm:approx}.

In the indicator-function setting considered here, \citet{doikov} start from the conceptual scheme that minimizes the true objective over the contracted domain and obtain Frank--Wolfe by replacing the objective with its first-order model. Higher-order models give contracting Newton and tensor methods. Their choice of the contracted domain is motivated by affine invariance, which the Euclidean ball lacks.

The derivation here makes the same two approximations in the opposite order: it linearizes first, inside the Euclidean ball of \algname{BPM}, which gives \algname{BroxGD}, and contracts second. It arrives at the same region as the smallest compact convex region inside the ball intersection, fixed independently of the gradient, that preserves \eqref{eq:comparison} (Proposition~\ref{prop:minimality}); Sections~\ref{sec:discover} and~\ref{sec:geometry} and Appendix~\ref{sec:geometric-appendix} measure what the approximation loses relative to the exact local step and how refinement can restore the distance guarantee.
\begingroup\small\setlength{\bibsep}{5pt}
\bibliographystyle{plainnat}
\bibliography{references}

@book{bertsekas,
  author    = {Dimitri P. Bertsekas},
  title     = {Nonlinear Programming},
  publisher = {Athena Scientific},
  address   = {Belmont, MA},
  edition   = {second},
  year      = {1999}
}

@book{boyd,
  author    = {Stephen Boyd and Lieven Vandenberghe},
  title     = {Convex Optimization},
  publisher = {Cambridge University Press},
  address   = {Cambridge},
  year      = {2004}
}

@article{doikov,
  author  = {Nikita Doikov and Yurii Nesterov},
  title   = {Affine-invariant contracting-point methods for convex optimization},
  journal = {Mathematical Programming},
  volume  = {198},
  number  = {1},
  pages   = {115--137},
  year    = {2023}
}

@article{dunn,
  author  = {Joseph C. Dunn and S. Harshbarger},
  title   = {Conditional gradient algorithms with open loop step size rules},
  journal = {Journal of Mathematical Analysis and Applications},
  volume  = {62},
  number  = {2},
  pages   = {432--444},
  year    = {1978}
}

@techreport{ferris,
  author      = {Michael C. Ferris and Sergei K. Zavriev},
  title       = {The linear convergence of a successive linear programming algorithm},
  institution = {Computer Sciences Department, University of Wisconsin--Madison},
  number      = {96-12},
  year        = {1996},
  url         = {https://pages.cs.wisc.edu/~ferris/techreports/96-12.pdf}
}

@article{frankwolfe,
  author  = {Marguerite Frank and Philip Wolfe},
  title   = {An algorithm for quadratic programming},
  journal = {Naval Research Logistics Quarterly},
  volume  = {3},
  number  = {1--2},
  pages   = {95--110},
  year    = {1956}
}

@article{garber,
  author  = {Dan Garber and Elad Hazan},
  title   = {A linearly convergent variant of the conditional gradient algorithm under strong convexity, with applications to online and stochastic optimization},
  journal = {SIAM Journal on Optimization},
  volume  = {26},
  number  = {3},
  pages   = {1493--1528},
  year    = {2016}
}

@article{brox,
  author  = {Kaja Gruntkowska and Hanmin Li and Aadi Rane and Peter Richt{\'a}rik},
  title   = {The ball-proximal (=``broximal'') point method: a new algorithm, convergence theory, and applications},
  journal = {arXiv preprint arXiv:2502.02002},
  year    = {2025},
  url     = {https://arxiv.org/abs/2502.02002}
}

@article{guelat,
  author  = {Jacques Gu{\'e}lat and Patrice Marcotte},
  title   = {Some comments on {W}olfe's `away step'},
  journal = {Mathematical Programming},
  volume  = {35},
  number  = {1},
  pages   = {110--119},
  year    = {1986}
}

@inproceedings{jaggi,
  author    = {Martin Jaggi},
  title     = {Revisiting {F}rank--{W}olfe: Projection-free sparse convex optimization},
  booktitle = {Proceedings of the 30th International Conference on Machine Learning},
  series    = {Proceedings of Machine Learning Research},
  volume    = {28},
  pages     = {427--435},
  year      = {2013},
  url       = {https://proceedings.mlr.press/v28/jaggi13.html}
}

@article{lacoste,
  author  = {Simon Lacoste-Julien and Martin Jaggi},
  title   = {An affine invariant linear convergence analysis for {F}rank--{W}olfe algorithms},
  journal = {arXiv preprint arXiv:1312.7864},
  year    = {2013},
  url     = {https://arxiv.org/abs/1312.7864}
}

@book{nemirovski,
  author    = {Arkadi Nemirovski and David Yudin},
  title     = {Problem Complexity and Method Efficiency in Optimization},
  publisher = {Wiley},
  address   = {New York},
  year      = {1983}
}

@article{nesterov2018,
  author  = {Yurii Nesterov},
  title   = {Complexity bounds for primal-dual methods minimizing the model of objective function},
  journal = {Mathematical Programming},
  volume  = {171},
  number  = {1--2},
  pages   = {311--330},
  year    = {2018}
}

@misc{warsawtalk,
  author       = {Peter Richt{\'a}rik},
  title        = {How to think about optimization: The broximal point method and its impact on optimization theory and machine learning},
  howpublished = {ML in PL, Warsaw, October 10, 2026. Talk slides},
  year         = {2026}
}

@article{sharper,
  author  = {Peter Richt{\'a}rik and Hanmin Li},
  title   = {A sharper theory of ball-proximal optimization: Convergence and radius selection},
  journal = {arXiv preprint arXiv:2609.35147},
  year    = {2026},
  url     = {https://arxiv.org/abs/2609.35147}
}

@misc{broxgd,
  author = {Peter Richt{\'a}rik and Kaja Gruntkowska and Hanmin Li},
  title  = {Broximal gradient descent: A projection-free sister of projected gradient descent},
  year   = {2026},
  url    = {https://arxiv.org/abs/2605.08850},
  note   = {Author manuscript, version 1.50; revision of arXiv:2605.08850}
}

@book{rockafellar,
  author    = {R. Tyrrell Rockafellar},
  title     = {Convex Analysis},
  publisher = {Princeton University Press},
  address   = {Princeton, NJ},
  year      = {1970}
}
\endgroup
\clearpage
\appendix
\numberwithin{theorem}{section}
\numberwithin{equation}{section}
\section{Residual bounds and refinement of the local solve}\label{sec:geometric-appendix}
Appendix~\ref{sec:geometric-appendix} supplies the definitions and proofs behind the outline in Section~\ref{sec:extension}. Appendix~\ref{sec:optimality} derives the optimality conditions of the exact local step and defines the residual. Appendix~\ref{sec:residual} proves the three-point inequality under a residual bound and the constant-radius guarantee that follows from it. Appendix~\ref{sec:refine} describes the refinement procedure on polytopes. Throughout, $\X$ is nonempty, closed, and convex, $x\in\X$, $g\in\R^n$, and $t>0$. Compactness of $\X$ is assumed only where stated.

\subsection{Optimality conditions of the exact local step}\label{sec:optimality}
Write the local problem \eqref{eq:brox} as
\begin{equation}\label{eq:local-program}
 \min_{z\in\X}F(z)\quad\text{subject to}\quad c(z)\le0,\qquad
 F(z)\eqdef\ip g{z-x},\quad c(z)\eqdef\tfrac12\bigl(\norm{z-x}^2-t^2\bigr).
\end{equation}
Only the ball constraint is dualized; the set $\X$ remains the domain over which the Lagrangian
\begin{equation}\label{eq:lagrangian}
 \mathcal L(z,\lambda)\eqdef F(z)+\lambda c(z),\qquad\lambda\ge0,
\end{equation}
is minimized. This keeps the global oracle available inside the analysis and avoids an explicit inequality description of $\X$.

\begin{lemma}[Optimality conditions of the local problem]\label{lem:kkt}
Let $\X$ be nonempty, closed, and convex, $x\in\X$, $t>0$, and $g\in\R^n$.
\begin{enumerate}
\item The problem \eqref{eq:local-program} has a solution $b$, and $\ip g{x-b}=M_t(x)$.
\item There is a multiplier $\lambda\ge0$ such that $b$ minimizes $\mathcal L(\cdot,\lambda)$ over $\X$ and $\lambda\,c(b)=0$. In particular $\lambda=0$ if $\norm{b-x}<t$.
\item For such $\lambda$,
\begin{equation}\label{eq:kkt-inequality}
 \ip{g+\lambda(b-x)}{z-b}\ge0\qquad(z\in\X).
\end{equation}
\end{enumerate}
Conversely, if $b\in\X\cap\B(x,t)$ and $\lambda\ge0$ satisfy \eqref{eq:kkt-inequality} and $\lambda\,c(b)=0$, then $b$ solves \eqref{eq:local-program}.
\end{lemma}
\begin{proof}
(1) The feasible set $\X\cap\B(x,t)$ contains $x$, is closed as an intersection of closed sets, and is bounded, hence compact. The objective is continuous, so a minimizer exists, and $\ip g{x-b}=M_t(x)$ by \eqref{eq:local-value}.

(2) Define the dual function $q(\lambda)\eqdef\inf_{z\in\X}\mathcal L(z,\lambda)$ for $\lambda\ge0$. For feasible $z$ and $\lambda\ge0$, $\lambda c(z)\le0$, so $q(\lambda)\le\mathcal L(z,\lambda)\le F(z)$: every dual value is a lower bound on $F(b)$. To establish strong duality and multiplier attainment, we verify Slater's condition by finding a point in the relative interior of $\X$ strictly inside the ball. Here $\operatorname{ri}\X$ denotes the interior of $\X$ within its affine hull, the smallest affine subspace containing $\X$. Although $x$ lies strictly inside its own ball, it may lie on the relative boundary of $\X$, so we construct a nearby relative-interior point. A nonempty convex set in $\R^n$ has a nonempty relative interior \citep[Theorem~6.2]{rockafellar}, so choose $\xi\in\operatorname{ri}\X$. If $x\in\operatorname{ri}\X$, take $\widehat z=x$. Otherwise take $\widehat z=(1-\alpha)x+\alpha\xi$ with $\alpha\in(0,1)$; then $\widehat z\in\operatorname{ri}\X$ \citep[Theorem~6.1]{rockafellar}, and $\norm{\widehat z-x}=\alpha\norm{\xi-x}<t$ for $\alpha$ small. Since $c$ is convex and $c(\widehat z)<0$, the refined Slater condition holds for the convex problem \eqref{eq:local-program} with domain $\X$, so the dual optimum is attained and equals the primal optimum \citep[Section~5.2.3]{boyd}: there is $\lambda\ge0$ with $q(\lambda)=F(b)$. Then
\[
 q(\lambda)\le\mathcal L(b,\lambda)=F(b)+\lambda c(b)\le F(b)=q(\lambda),
\]
so both inequalities are equalities: $b$ minimizes $\mathcal L(\cdot,\lambda)$ over $\X$, and $\lambda c(b)=0$. If $\norm{b-x}<t$, then $c(b)<0$ and $\lambda=0$.

(3) The function $\mathcal L(\cdot,\lambda)$ is convex and differentiable, with $\nabla_z\mathcal L(b,\lambda)=g+\lambda(b-x)$. Fix $z\in\X$. Convexity of $\X$ puts $b+\alpha(z-b)$ in $\X$ for $\alpha\in(0,1]$, and minimality of $b$ gives
\[
 0\le\frac{\mathcal L(b+\alpha(z-b),\lambda)-\mathcal L(b,\lambda)}{\alpha}
 \ \xrightarrow[\alpha\downarrow0]{}\ \ip{\nabla_z\mathcal L(b,\lambda)}{z-b},
\]
which is \eqref{eq:kkt-inequality}.

For the converse, let $z\in\X\cap\B(x,t)$. Since $\lambda c(z)\le0$, convexity of $\mathcal L(\cdot,\lambda)$ gives
\[
 F(z)\ge\mathcal L(z,\lambda)\ge\mathcal L(b,\lambda)+\ip{g+\lambda(b-x)}{z-b}
 \overset{\eqref{eq:kkt-inequality}}{\ge}\mathcal L(b,\lambda)=F(b)+\lambda c(b)=F(b).
\]
\end{proof}

The inequality \eqref{eq:kkt-inequality} is the property of the exact step that the distance argument uses, and its violation by a candidate can be measured. For a candidate $y\in\X$ and $\lambda\ge0$, the residual
\begin{equation}\label{eq:residual}
 \varepsilon_\lambda(x,y)\eqdef\sup_{z\in\X}\ip{g+\lambda(y-x)}{y-z}
\end{equation}
is the largest first-order decrease of $\mathcal L(\cdot,\lambda)$ available from $y$ over $\X$, since $\nabla_z\mathcal L(y,\lambda)=g+\lambda(y-x)$. It is nonnegative, since $z=y$ is admissible, and it vanishes exactly when $y$ satisfies \eqref{eq:kkt-inequality} in place of $b$, which for the convex function $\mathcal L(\cdot,\lambda)$ means that $y$ minimizes it over $\X$.

When the linear minimum is attained, one global oracle call in that direction computes the residual: if $v$ minimizes $\ip{g+\lambda(y-x)}{z}$ over $\X$, then $\varepsilon_\lambda(x,y)=\ip{g+\lambda(y-x)}{y-v}$. By Lemma~\ref{lem:kkt}(3), the exact step with its multiplier has residual zero. The residual depends on the choice of $\lambda$, even at the exact step: on $\X=[-2,2]$ with $x=0$, $g=-1$, and $t=1$, the exact step $b=1$ has residual zero for its multiplier $\lambda=1$, but residual one for $\lambda=0$.

\subsection{From the residual to a constant-radius guarantee}\label{sec:residual}
We next show how a nonzero residual changes the three-point inequality used in the distance analysis of \algname{BroxGD}.

\begin{theorem}[Inexact three-point inequality]\label{thm:threepoint}
Let $\X$ be nonempty, closed, and convex, let $x,y\in\X$, and let $g\in\R^n$. Suppose $\varepsilon_\lambda(x,y)\le\varepsilon<\infty$ for some $\lambda\ge0$. Then, for every $z\in\X$,
\begin{equation}\label{eq:threepoint}
 \ip g{y-z}\le\varepsilon+\frac\lambda2
 \bigl(\norm{x-z}^2-\norm{y-z}^2-\norm{y-x}^2\bigr).
\end{equation}
\end{theorem}
\begin{proof}
Set $d\eqdef y-x$. The residual bound gives, for every $z\in\X$,
\begin{equation}\label{eq:residual-comparison}
 \ip g{y-z}\overset{\eqref{eq:residual}}{\le}\varepsilon-\lambda\ip d{y-z}.
\end{equation}
The identity
\begin{equation}\label{eq:threepoint-identity}
 2\ip d{y-z}=\norm{y-z}^2-\norm{x-z}^2+\norm d^2
\end{equation}
follows by expanding $\norm{x-z}^2=\norm{(y-z)-d}^2$. Substituting it into \eqref{eq:residual-comparison} gives \eqref{eq:threepoint}.

\end{proof}

For the exact local step with its multiplier, $\varepsilon=0$ by Lemma~\ref{lem:kkt}, and \eqref{eq:threepoint} is the three-point inequality used in the proof of \citet[Theorem~3.1]{broxgd}.

The residual also bounds the local-model error. If $y\in\mathcal C_t(x)$ and $\varepsilon_\lambda(x,y)<\infty$ for some $\lambda\ge0$, applying \eqref{eq:threepoint} at the exact local step $b$ gives
\begin{equation}\label{eq:cert-local-error}
 0\le M_t(x)-\ip g{x-y}
 =\ip g{y-b}
 \le\varepsilon_\lambda(x,y)+\frac\lambda2\bigl(t^2-\norm{y-x}^2\bigr).
\end{equation}
Here we used $\norm{b-x}\le t$ and discarded the nonpositive term $-\lambda\norm{y-b}^2/2$. For a full-radius candidate, the local-model error is at most the residual itself.

We now use a residual test to recover the constant-radius objective bound of \algname{BroxGD} without a dependence on the diameter.

\begin{corollary}[Constant radius with a residual test]\label{cor:certified}
Let $\X$ be nonempty, closed, and convex. Let $f$ satisfy Assumption~\ref{ass:objective}, and fix a minimizer $x_\star$ and an initial point $x_0\in\X$. Choose $0\le\theta<1$, $R\ge\norm{x_0-x_\star}$ with $R>0$, and an integer $K\ge1$. Use the constant radius
\begin{equation}\label{eq:cert-radius}
 t\eqdef\frac{R}{\sqrt{(1-\theta)K}}.
\end{equation}
At each step, with $g=\nabla f(x)$, return a feasible $y$ satisfying either
\begin{equation}\label{eq:cert-test}
 \norm{y-x}=t,\qquad\lambda>0,\qquad
 \varepsilon_\lambda(x,y)\le\frac{\theta\lambda t^2}{2},
\end{equation}
or $\norm{y-x}\le t$ and $\varepsilon_0(x,y)=0$, in which case stop. Otherwise set $x\leftarrow y$ and continue until $K$ steps have been performed. The best objective gap among the initial point and the returned points satisfies
\begin{equation}\label{eq:cert-rate}
 \min\Delta\le\frac{LR^2}{2(1-\theta)K}.
\end{equation}
\end{corollary}
\begin{proof}
Smoothness \eqref{eq:smooth} at $y$ and convexity \eqref{eq:convex} at $x_\star$, both based at $x$, give
\begin{equation}\label{eq:endpoint-gap}
 \Delta(y)\le\ip g{y-x_\star}+\frac L2\norm{y-x}^2.
\end{equation}
Set $\sigma\eqdef Lt^2/2$. We show that either a returned point has objective gap at most $\sigma$, or each step decreases the squared distance to $x_\star$ by more than $(1-\theta)t^2$.

If the run stops with $\lambda=0$ and $\varepsilon_0(x,y)=0$, then \eqref{eq:endpoint-gap} and \eqref{eq:residual} with $\lambda=0$ give $\ip g{y-x_\star}\le0$, so $\Delta(y)\le\frac L2\norm{y-x}^2\le\sigma$.

Otherwise $y$ satisfies \eqref{eq:cert-test}. If $\Delta(y)>\sigma$, then \eqref{eq:endpoint-gap} with $\norm{y-x}=t$ gives $\ip g{y-x_\star}\ge\Delta(y)-\sigma>0$. Rearranging \eqref{eq:threepoint} at $z=x_\star$, with $\norm{y-x}=t$ and $\varepsilon=\varepsilon_\lambda(x,y)$,
\[
 \norm{y-x_\star}^2
 \le\norm{x-x_\star}^2-t^2+\frac{2\varepsilon_\lambda(x,y)}{\lambda}-\frac{2\ip g{y-x_\star}}{\lambda}
 \overset{\eqref{eq:cert-test}}{<}\norm{x-x_\star}^2-(1-\theta)t^2.
\]
If all $K$ returned points had gap greater than $\sigma$, summing these strict decreases would give
\[
 \norm{x_K-x_\star}^2<R^2-K(1-\theta)t^2\overset{\eqref{eq:cert-radius}}{=}0,
\]
which is impossible. So some returned point has gap at most $\sigma$, and \eqref{eq:cert-radius} gives \eqref{eq:cert-rate}.
\end{proof}

The residual test requires neither $x_\star$ nor $f_\star$. The bound $R$ enters only the radius. The condition $\Delta(y)>\sigma$ in the proof plays the role of the Type-I radius condition \eqref{eq:type-one}. Under \eqref{eq:type-one}, every locally feasible $y$ satisfies
\[
 \ip g{y-x_\star}\ge\ip g{x-x_\star}-\norm g\,t\ge0.
\]
For a candidate that also satisfies the residual test \eqref{eq:cert-test}, the same argument gives
\[
 \norm{y-x_\star}^2\le\norm{x-x_\star}^2-(1-\theta)t^2.
\]
The corollary uses the objective threshold instead, because \eqref{eq:type-one} involves the unknown $x_\star$.

\subsection{Refinement on a polytope}\label{sec:refine}
Let $\X=\operatorname{conv}V$ with $V$ a finite set of vertices, and let the global oracle return a minimizing vertex. Fix $x$, $g$, $t>0$, and $\theta\in[0,1)$. The following procedure produces a point passing the test \eqref{eq:cert-test}, or a point with $\lambda=0$ and zero residual.
\begin{enumerate}
\item Call the oracle with direction $g$ to obtain $s$. If $s=x$, return $x$; its residual with $\lambda=0$ is zero. Otherwise set $P_0\eqdef\operatorname{conv}\{x,s\}$ and $j\eqdef0$.
\item Solve the restricted problem
\[
 \min_{y\in P_j\cap\B(x,t)}\ip g{y-x}
\]
and obtain a multiplier $\lambda\ge0$ satisfying
\[
 \begin{aligned}
 \ip{g+\lambda(y-x)}{z-y}&\ge0\qquad(z\in P_j),\\
 \lambda\bigl(\norm{y-x}^2-t^2\bigr)&=0.
 \end{aligned}
\]
Both a solution and such a multiplier exist by Lemma~\ref{lem:kkt} applied to $P_j$.
\item Call the oracle with direction $w\eqdef g+\lambda(y-x)$ to obtain a vertex $v$, and compute $\varepsilon\eqdef\ip w{y-v}=\varepsilon_\lambda(x,y)$.
\item If \eqref{eq:cert-test} holds, or if $\lambda=0$ and $\varepsilon=0$, return $y$. Otherwise set $P_{j+1}\eqdef\operatorname{conv}(P_j\cup\{v\})$, increase $j$, and go to step 2.
\end{enumerate}
The first restricted solution can be taken as
\begin{equation}\label{eq:first-ray}
 y_0\eqdef x+\min\left\{1,\frac{t}{\norm{s-x}}\right\}(s-x),
\end{equation}
the truncated Frank--Wolfe step along the oracle direction. At $t=\gamma D$ it may move farther than the contracted step, but its linear value is at least as small. The sets $P_j$ grow with $j$ and stay inside $\X$, so the restricted optimal values cannot increase. Refinement improves the linear model at a fixed radius, but need not improve the actual objective relative to the first candidate.

\begin{proposition}[Finite refinement preserves the contracted comparison]\label{prop:refine}
The procedure terminates after at most $|V|+1$ global oracle calls, including initialization and the final check. This bound does not include the cost of the restricted solves. Every restricted candidate has step length at most $t$ and linear objective no greater than that of $y_0$. If $D>0$ bounds the diameter of $\X$ and $t=\gamma D$ with $0<\gamma\le1$, every candidate satisfies \eqref{eq:approx-condition} with $e_k=0$ when taken as the next iterate. Under Assumption~\ref{ass:objective}, with $g=\nabla f(x)$, each candidate consequently satisfies the recurrence \eqref{eq:exact-recurrence}.
\end{proposition}
\begin{proof}
Termination. Suppose the test fails at step 4. If $\lambda>0$, complementary slackness gives $\norm{y-x}=t$, so failure of \eqref{eq:cert-test} means $\varepsilon>\theta\lambda t^2/2\ge0$; if $\lambda=0$, failure means $\varepsilon>0$. In both cases $\ip w{y-v}>0$. Restricted optimality gives $\ip w{y-z}\le0$ for every $z\in P_j$, so $v\notin P_j$: every failed test adds a vertex of $V$ that was absent. Initialization supplies the vertex $s$, at most $|V|-1$ failed tests add the remaining vertices, and one final check remains. Once $P_j=\X$, the restricted solution is the exact local step and its residual is zero by Lemma~\ref{lem:kkt}, so the test passes.

Step length and linear value. Every candidate lies in $\B(x,t)$, and $P_j\supseteq P_0$ gives $\ip g{y-x}\le\ip g{y_0-x}$.

Comparison. With $G=\ip g{x-s}\ge0$ and $\norm{s-x}\le D$, the fraction in \eqref{eq:first-ray} is at least $\gamma$, so for every $z\in\X$,
\[
 \ip g{y-x}\le\ip g{y_0-x}\overset{\eqref{eq:first-ray}}{\le}-\gamma G
 \overset{\eqref{eq:oracle}}{\le}\gamma\ip g{z-x}.
\]
This is \eqref{eq:approx-condition} with zero error, and Theorem~\ref{thm:approx} gives the recurrence when $g=\nabla f(x)$.
\end{proof}

For the triangle family of Section~\ref{sec:geometry}, take any $h>0$ and $0<t\le1$.
The initial oracle returns $s=(2,h)$, and the first restricted solution and its multiplier are
\[
 y_0=\frac{t}{D_h}(2,h),\qquad
 \lambda_0=\frac{2}{tD_h}.
\]
The residual direction is
\[
 w=g+\lambda_0y_0=\frac1{D_h^2}(-h^2,2h).
\]
Its inner products with $(0,0)$ and $(2,h)$ are zero, while $\ip w{(1,0)}=-h^2/D_h^2$.
Thus the oracle call selects $v=(1,0)$, and the residual and test threshold are
\[
 \varepsilon=\ip w{y_0-v}=\frac{h^2}{D_h^2},\qquad
 \frac{\theta\lambda_0t^2}{2}=\frac{\theta t}{D_h}.
\]
If the test fails, adding $v$ gives $P_1=\X_h$.
The next restricted solve returns the exact \algname{BroxGD} step $b=(t,0)$ with multiplier $\lambda=1/t$, direction $w=0$, and residual zero.
Hence one refinement suffices whenever the initial candidate fails the test.

In particular, every distance-increasing step identified in Section~\ref{sec:geometry} fails the test.
Indeed, $D_h>4/t\ge4$ implies
\[
 h^2=D_h^2-4>D_h>\theta tD_h,
\]
so $\varepsilon>\theta t/D_h$ for every $\theta\in[0,1)$.
The refinement therefore recovers the distance guarantee throughout this family of counterexamples.
At $t=1$ it reaches $x_\star$; for smaller radii it returns $b=(t,0)$.

\section{A nuclear-norm example}\label{sec:nuclear-example}
This example supports the computational comparison in the introduction. For matrices, use the Frobenius inner product $\ip G Z=\operatorname{tr}(G^\top Z)$ and norm $\norm{Z}_{\mathrm F}$; the nuclear norm $\norm{Z}_*$ is the sum of the singular values. Take the current point $x=0$ and the feasible set
\[
 \X=\{Z\in\R^{r\times r}:\norm{Z}_*\le1\},\qquad r\ge2.
\]
Choose $a_1>\cdots>a_r>0$ with $\sum_{i=1}^r a_i=1$, and orthogonal matrices $U,V$. Set
\[
 A=U\operatorname{diag}(a_1,\ldots,a_r)V^\top,\qquad
 W=UV^\top,\qquad G=W+A,\qquad
 t=\norm{A}_{\mathrm F}=\Bigl(\sum_{i=1}^r a_i^2\Bigr)^{1/2}<1.
\]
The singular values of $G$ are $\sigma_i=1+a_i$. Thus the global oracle returns $s=-u_1v_1^\top$, using only a leading singular-vector pair \citep[Table~1]{jaggi}. The Frobenius diameter of $\X$ is $D=2$, since $\norm{Z}_{\mathrm F}\le\norm{Z}_*\le1$ on $\X$ and $\pm u_1v_1^\top$ attain it. With $x=0$ and $\gamma=t/2$, the Frank--Wolfe step at radius $t$ is
\[
 y=-\frac t2u_1v_1^\top.
\]
Linear minimization over the Frobenius ball alone returns $-tG/\norm{G}_{\mathrm F}$. The local problem is
\begin{equation}\label{eq:nuclear-local}
 \min_Z\{\ip G Z:\norm{Z}_*\le1,\ \norm{Z}_{\mathrm F}\le t\}.
\end{equation}
For every feasible $Z$, nuclear--operator norm duality and Cauchy--Schwarz give
\[
 \begin{aligned}
 -\ip G Z
 &=-\ip W Z-\ip A Z\\
 &\le\norm{W}_{\mathrm{op}}\norm{Z}_*
       +\norm{A}_{\mathrm F}\norm{Z}_{\mathrm F}
 \le1+t^2,
 \end{aligned}
\]
where $\norm{W}_{\mathrm{op}}=1$. The point $b=-A$ attains this bound and satisfies both constraints with equality, so it solves \eqref{eq:nuclear-local}. Equality throughout forces $Z=-A$, so the solution is unique. In terms of the singular values of $G$,
\[
 b=-U\operatorname{diag}\bigl((\sigma_1-1)_+,\ldots,(\sigma_r-1)_+\bigr)V^\top,
 \qquad (c)_+\eqdef\max\{c,0\}.
\]
All $r$ coefficients are positive, so the \algname{BroxGD} step $b$ has rank $r$, whereas the Frank--Wolfe step $y$ has rank one. The direct spectral computation of $b$ uses every singular component of $G$, whereas the global oracle uses only the leading pair.

\end{document}